\documentclass[11pt]{article}
\usepackage{a4}
\usepackage{latexsym}
\usepackage{amsfonts}
\usepackage{amssymb}
\usepackage{amscd}
\usepackage{amsmath}
\usepackage{epsfig}
\begin{document}
\title{\textbf{Poincar\'{e} series of a toric variety}}
\author{Ann Lemahieu \date{}\footnote{Partially supported by the grants FWO G.0318.06, FWO V
4.034.06N and MEC PN I+D+I MTM2004-00958}} \maketitle {\footnotesize
\emph{\textbf{Abstract.---} For an affine toric variety $X$ we
compute the Poincar\'{e} series of the multi-index filtration
defined by a finite number of monomial divisorial valuations on the
ring $\mathcal{O}_{X,0}$. We give an alternative description of the
Poincar\'{e} series as an integral with respect to the Euler
characteristic over the projectivization of the space of germs
$\mathcal{O}_{X,0}$.
\\In particular we study divisorial valuations on the ring
$\mathcal{O}_{\mathbb{C}^d,0}$ that arise by considering toric
constellations. We give an explicit formula for the Poincar\'{e}
series and a nice geometric description. This generalizes an
expression of the Poincar\'{e} series for curves and rational
surface singularities.}}
\\  ${}$
\\ \\ \enlargethispage{1cm}
\indent A. Campillo, F. Delgado and S.M. Gusein-Zade computed the
Poincar\'{e} series of the multi-index filtration on the ring of
germs of functions in two variables defined by orders of a function
on the branches of a reducible curve singularity in \cite{trio2} and
\cite{trio5}. For plane curves they showed that the Poincar\'{e}
series coincides with the Alexander polynomial of the link of the
singularity. An intuitive motivation for this phenomenon is still
missing. Also for an arbitrary collection of plane divisorial
valuations the Poincar\'{e} series has been studied, see \cite{duo}.
In \cite{trio6} the Poincar\'{e} series of the multi-index
filtration on the ring of germs of functions on a rational surface
singularity defined by the multiplicities of a function along
components of the
exceptional divisor of a resolution of the singularity has been computed. \\
\indent These Poincar\'{e} series can be written in several ways.
They can be described by the fibers of the corresponding extended
semigroups. They also have the shape of an integral with respect to
the Euler characteristic over the projectivization of the space of
germs of functions (see \cite{trio4} and \cite{trio6}). This notion
is similar to the notion of motivic integration and has been
introduced in \cite{trio3}. Furthermore, one has a description at
the level of the modification space.
\\ \\ \indent As the following shows, for an arbitrary affine variety $X$ with a finite set of divisorial
valuations on the ring of germs $\mathcal{O}_{X,o}$, with $o$ a
point on $X$, the multi-index filtration induced by the valuations
does not always permit to define its Poincar\'{e} series.
\\ Let $k$ be an arbitrary field. The Poincar\'{e} series of $X$ is
given by a multi-index filtration with index set $\mathbb{Z}^r$ that
is induced by a set of discrete valuations $\{\nu_1,\cdots,\nu_r\}$
of $\mathcal{O}_{X,o}$. For $\underline{v}=(v_1,\cdots,v_r) \in
\mathbb{Z}^r$, such valuations induce complete ideals
$I(\underline{v}):= \{g \in \mathcal{O}_{X,o} \mid \nu_i(g) \geq
v_i, 1 \leq i \leq r\}$. When the centre of each valuation $\nu_i$,
i.e. the set $\{g \in \mathcal{O}_{X,o} \mid \nu_i(g) > 0 \}$, is
the maximal ideal $\mathsf{m}$ of $\mathcal{O}_{X,o}$ ($1 \leq i
\leq r$), the $k$-vector spaces
$I(\underline{v})/I(\underline{v}+\underline{1})$ have finite
dimension. Let us denote
$d(\underline{v}):=\mbox{dim}I(\underline{v})/I(\underline{v}+\underline{1})$.
Then the series $L(t_1,\cdots,t_r)=\sum_{\underline{v} \in
\mathbb{Z}^r} d(\underline{v}) \underline{t}^{\underline{v}}$ is a
well-defined Laurent series. Let $\underline{v}$ and
$\underline{v}'$ be vectors in $\mathbb{Z}^r$, let $j \in
J:=\{1,\cdots, r\}$ and suppose that $v_i=v_i'$ for all $i \in
J\setminus j$. If $v_j \leq 0$ and $v_j' \leq 0$, then
$I(\underline{v})=I(\underline{v}')$ and hence
$\prod_{i=1}^{r}(t_i-1)L(t_1,\cdots,t_r) \in
\mathbb{Z}[[t_1,\cdots,t_r]]$. The Poincar\'{e} series is then
defined as
\[P(t_1,\cdots,t_r)=\frac{\prod_{i=1}^{r}(t_i-1)L(t_1,\cdots,t_r)}{(t_1\cdots
t_r - 1)}.\] However, if one of the valuations $\nu_1,\cdots,\nu_r$
does not have its centre at $\mathsf{m}$, then one can not define
the series $L$. Indeed, suppose $\nu_i$ ($1 \leq i \leq r$) is a
valuation with centre at the prime ideal $\mathsf{p}_i$ which is
different from $\mathsf{m}$. If $g \in \mathcal{O}_{X,o}$ and
$\underline{\nu}(g)=\underline{v}$, then choose a function $h \in
\mathsf{m} \setminus \mathsf{p}_i$. Now for each $n \in
\mathbb{Z}_{\geq 0}$ one has that $gh^n \in I(\underline{v})$ and
$gh^n \notin I(\underline{v}+\underline{1})$ and, as all the $gh^n
(n \in \mathbb{Z}_{\geq 0})$ are linearly independent over $k$,
it follows that $d(\underline{v})$ is infinite. 
\\ \\
\indent In this article we work over the field $k=\mathbb{C}$. We
compute the Poincar\'{e} series for affine toric varieties and
provide some equivalent descriptions for the Poincar\'{e} series
such as exist in the cases of curves, rational surface singularities
and plane divisorial valuations. As an interesting example, we study
the case of divisorial valuations on the ring
$\mathcal{O}_{\mathbb{C}^d,o} (d\geq 2)$ that are created by toric
constellations.
\\ \\  ${}$ \begin{center} \textsc{1. Poincar\'{e} series of an affine toric variety}
\end{center} ${}$\\
\indent 
Let $\sigma \subset N \otimes_{\mathbb{Z}}\mathbb{R}$ be a rational
finite polyhedral strongly convex $d$-dimensional cone, where
$N=\mathbb{Z}^d$. Let $M$ be the dual space to $N$, then there is a
natural bilinear map
$M \times N  \rightarrow  \mathbb{Z}: (m,n)  \mapsto \langle m,n
\rangle$. The dual cone $\check{\sigma}$ to $\sigma$ is defined as
the set $\{m \in M \mid \langle m,x \rangle \geq 0, \forall x \in
\sigma\}$. When $S$ is a semigroup in $\check{\sigma} \cap M$ that
generates $\check{\sigma} \cap M$ as cone, then also $S$ is finitely
generated. Without loss of generality, we assume that $M$ is also
the group generated by the semigroup $S$. We consider
$\mathbb{C}[S]$ as the $S$-graded algebra $\oplus_{s \in
S}\mathbb{C} \chi^s$. Let $X$ be the affine toric variety
Spec$(\mathbb{C}[S])$. Note that $X$ is a normal variety if and only
if $S=\check{\sigma} \cap M$. Let $\pi:X' \rightarrow X$ be an
equi\-variant proper birational morphism, $X'$ being another toric
variety. An irreducible codimension $1$ subvariety $D$ of $X'$
induces a discrete valuation $\nu_D$ of $\mathbb{C}(X)$ and so an
element $n_D$ of $N$ that acts as follows:
\begin{eqnarray*} n_D: M
& \longrightarrow & \mathbb{Z}\\
m & \longmapsto & \nu_D(\chi^m).
\end{eqnarray*}
If $D$ is an irreducible component of $X' \setminus T$ then $n_D$ is
a primitive element in $\sigma \cap N$. Vice versa a primitive
element $n$ in $\sigma \cap N$ defines a discrete valuation $\nu$ of
$\mathbb{C}(X)$ by setting $\nu(\sum_{m \in F} a_m\chi^m)=min\{n(m)
\mid m \in F, a_m \neq 0\}$. In fact, one has $\nu=\nu_D$ for some
$D$ in $X'$ for an appropriated $\pi: X' \rightarrow X$. Such
valuations are usually called monomial valuations. \\ As we saw
before, the valuations inducing the Poincar\'{e} series must have
the unique $0$-dimensional orbit of $X$ as centre. If $D$ is an
irreducible component of codimension $1$ of $X' \setminus T$, then
$\pi(D)$ is the closure of the orbit that is associated to the
unique face $\tau$ of $\sigma$ such that $\overset{\circ}{\tau}$
contains $n_D$. This means that the valuations that we are
considering correspond bijectively to the primitive elements of
$\overset{\circ}{\sigma} \cap N$. In what follows we identify both
notions.  The above correspondence also stands for non-primitive
elements in $\overset{\circ}{\sigma} \cap N$ and non-normalized
monomial valuations.
\\ \\ \indent We compute the Poincar\'{e} series $P_{X,H}$ of the affine toric variety $X$ with respect to an
arbitrary set of divisorial monomial valuations
$H=\{\nu_1,\cdots,\nu_r\}$ such that $H \subset
\overset{\circ}{\sigma}$. In what follows $\underline{\nu}$ stands
for the vector $(\nu_1,\cdots,\nu_r)$ and by $\langle
s,\underline{\nu}\rangle$ we mean the vector $(\langle s, \nu_1
\rangle,\langle s, \nu_2\rangle,\cdots,\langle s,\nu_r\rangle)$.
We define the monomial cone $C$ as the set of values that can be
obtained by valuating monomials, i.e. $C=\{\underline{\nu}(m) \mid
m \mbox{ monomial in } \mathbb{C}[S]\}$.
\\One has a $\mathbb{Z}$-linear map $\phi: M
\longrightarrow \mathbb{Z}^r$ given by $\phi(m)=\langle
m,\underline{\nu}\rangle$. This map induces the following map among
Laurent series groups \small{\begin{eqnarray*}
\Phi: \mathbb{Z}[[M]]=\mathbb{Z}[[u_1,\cdots,u_d,u_1^{-1},\cdots,u_d^{-1}]] & \longrightarrow & \mathbb{Z}[[\mathbb{Z}^r]]=\mathbb{Z}[[t_1,\cdots,t_r,t_1^{-1},\cdots,t_r^{-1}]]\\
\sum_{i}\lambda_i\underline{u}^{m} & \longmapsto &
\sum_{i}\lambda_i\underline{t}^{\langle m,\underline{\nu} \rangle}.
\end{eqnarray*}}
\normalsize \noindent Notice that for each commutative group
$\Gamma$ the Laurent series group $\mathbb{Z}[[\Gamma]]$ is in fact
a $\mathbb{Z}[\Gamma]$-module, and that the assignment of this
module
to $\Gamma$ is functorial. \\
In particular, for our given semigroup $S$, one has the multi-graded
Poincar\'{e} series of commutative algebra $Q(\underline{u})=\sum_{s
\in S}\underline{u}^s$ which is an element of $\mathbb{Z}[[S]]$ and
so can be interpreted as an element in $\mathbb{Z}[[M]]$. In fact,
this multi-graded Poincar\'{e} series is usually expressed as a
rational function having $Q(\underline{u})$ as power series
expansion. The Poincar\'{e} series in geometry
$P_{X,H}(\underline{t})$ is an element in
$\mathbb{Z}[[\mathbb{N}^r]]$ and so also an element of
$\mathbb{Z}[[\mathbb{Z}^r]]$. The following result shows the
relationship between both Poincar\'{e} series.
\\ \\
\textbf{Theorem 1.---} \emph{The Poincar\'{e} series $P_{X,H}$
defined by the multi-index filtration induced by $H$ and associated
to the multi-graduation of $S$, is the image under $\Phi$ of the
Poincar\'{e} series of commutative algebra of the semigroup $S$,
i.e.}
\[P_{X,H}(\underline{t})=\Phi(Q(\underline{u})).\]
${}$\\
\emph{Proof.} For a set $A \subset \{i_1,\cdots,i_s\}$, let
$\alpha_A$ be the function
\begin{eqnarray*}
\alpha_A:  \mathbb{Z}^s & \longrightarrow &\mathbb{Z}^s \\
\underline{v} & \longmapsto & \underline{v}'
\end{eqnarray*}
where $v'_i=v_i-1$ if $i \in A$ and $v'_i=v_i$ if $i \notin A$.\\
Then, for $\underline{v} \in \mathbb{Z}^r$, the coefficient of
$\underline{t}^{\underline{v}}$ in
$\prod_{i=1}^{r}(t_i-1)L(t_1,\cdots,t_r)$ is \[(-1)^r \sum_{A
\subset \{1,\cdots,r\}} (-1)^{\#A} \mbox{dim}
\frac{I(\alpha_A(\underline{v}))}{I(\alpha_A(\underline{v})+
\underline{1})}\] and the coefficient of
$\underline{t}^{\underline{v}}$ in $P_{X,H}(\underline{t})$ is
\begin{eqnarray*}
(-1)^{r+1}\sum_{A \subset \{1,\cdots,r\}} (-1)^{\#A} \mbox{dim}
\frac{I(\alpha_A(\underline{v}))}{I(\alpha_A(\underline{v})+
\underline{1})} &+ &\\ (-1)^{r+1}\sum_{A \subset \{1,\cdots,r\}}
(-1)^{\#A} \mbox{dim}
\frac{I(\alpha_A(\underline{v})-\underline{1})}{I(\alpha_A(\underline{v}))}&
+& \\ (-1)^{r+1}\sum_{A \subset \{1,\cdots,r\}} (-1)^{\#A}
\mbox{dim}\frac{I(\alpha_A(\underline{v})-\underline{2})}{I(\alpha_A(\underline{v})-
\underline{1})} &+ & \cdots
\end{eqnarray*}
This finite sum can be rewritten as
\begin{eqnarray*}
 (-1)^{r+1} \sum_{A \subset \{1,\cdots,r\}} (-1)^{\#A} \mbox{dim}
\frac{\mathbb{C}[S]}{I(\alpha_A(\underline{v})+ \underline{1})}.
\end{eqnarray*}
Every subset $A \subset \{1,\cdots,r\}$ can be written in a unique
way as $A_1 \times A_2$, with $A_1 \subset \{2,\cdots,r\}$ and
$A_2 \subset \{1\}$. We group the terms having the same component
$A_1$ and we get \small{\begin{eqnarray*} & & \sum_{A \subset
\{1,\cdots,r\}} (-1)^{\#A} \mbox{dim}
\frac{\mathbb{C}[S]}{I(\alpha_A(\underline{v})+ \underline{1})} \\
&=& \sum_{A \subset \{2,\cdots,r\}} (-1)^{\# A}\#\{\chi^s \mid
\langle s,\nu_1\rangle=v_1, \langle s,(\nu_2,\cdots,\nu_r)\rangle
\geq \alpha_A(v_2,\cdots,v_r)+\underline{1}\}.
\end{eqnarray*}}
\normalsize \noindent We go on simplifying in the same way, so now
we write every subset $A \subset \{2,\cdots,r\}$ as $A_1 \times
A_2$, with $A_1 \subset \{3,\cdots,r\}$ and $A_2 \subset \{2\}$ and
so on. At each step we group the terms having the same component
$A_1$ and we obtain \small
\begin{eqnarray*}
 & & \sum_{A \subset \{3,\cdots,r\}} (-1)^{\# A+1}\#\{\chi^s \mid
\langle s,\nu_1\rangle=v_1,\langle s,\nu_2\rangle=v_2, \\ & &
 \qquad \qquad \qquad \qquad \qquad\langle s,(\nu_3,\cdots,\nu_r)\rangle \geq
\alpha_A(v_3,\cdots,v_r)+\underline{1}\} \\
&= & \sum_{A \subset \{4,\cdots,r\}} (-1)^{\# A+2}\#\{\chi^s \mid
\langle s,\nu_1\rangle=v_1,\langle s,\nu_2\rangle=v_2, \langle
s,\nu_3\rangle=v_3,\\ & &
 \qquad \qquad \qquad \qquad \qquad\langle s,(\nu_4,\cdots,\nu_r)\rangle \geq
\alpha_A(v_4,\cdots,v_r)+\underline{1}\}\\
& \vdots & \\
&=& (-1)^{r-1}\#\{\chi^s \mid \langle
s,\underline{\nu}\rangle=\underline{v}\}.
\end{eqnarray*}
\normalsize Hence the Poincar\'{e} series $P_{X,H}(\underline{t})$
is
\[\sum_{\underline{v} \in \mathbb{Z}^r} \#\{\chi^s \mid \langle
s,\underline{\nu}\rangle=\underline{v}\}
\underline{t}^{\underline{v}}.\] \hfill $\blacksquare$  \newpage
\noindent In what follows we will write $N(\underline{v})$ when we
want to refer to the number $\# \{s \in S \mbox{ $|$ }\langle
s,\underline{\nu}\rangle=\underline{v}\}$. Since $\sigma$ is
$d$-dimensional and strongly convex, one has that $N(\underline{v})$
is finite for every $\underline{v}$.
\\ \\ \indent The semigroup of $X$
with respect to the set of valuations $H$ is
\[S_{H}=\{\underline{v} \in \mathbb{Z}^r \mid \underline{v}=\underline{\nu}(g) \mbox{ for some } g \in \mathcal{O}_{X,o} \}.\]
Let $g$ be a function in $I(\underline{v})$, say
$g=\sum_{i=1}^{s}\lambda_i m_i$ where the $m_i$ are monomials and
the $\lambda_i$ are complex numbers ($1 \leq i \leq s$). For $j \in
J$, we denote by $a_j(g)$ the part $\sum_{i \in K_j}\lambda_i m_i$,
where $K_j=\{i \mid \nu_j(m_i)=v_j\}$, and by
$\underline{a}(g)=(a_1(g),\cdots,a_r(g))$. The set
$\hat{S}_{H}=\{(\underline{\nu}(g),\underline{a}(g)) \mid g \in
\mathcal{O}_{X,o}\}$ is a semigroup with respect to the summation of
the components $\underline{\nu}$ and multiplication of the parts
$\underline{a}$ and is called the extended semigroup. This notion
showed up for the first time in \cite{trio1} where it was introduced
for plane curves. Later this notion has been extended in the study
of the Poincar\'{e} series of plane divisorial
valuations and rational surface singularities. \\
For $j \in J$, denote by $D_j(\underline{v})$ the complex vector
space $I(\underline{v})/I(\underline{v}+\underline{e_j})$ where
$e_j$ is the $r$-tuple with \emph{j}-th component equal to $1$ and
the other components equal to $0$. Consider the application
\begin{eqnarray*}
j_{\underline{v}}: I(\underline{v}) & \longrightarrow &
D_1(\underline{v}) \times \cdots \times D_r(\underline{v}) \\
g & \longmapsto & (a_1(g),\cdots,a_r(g))=\underline{a}(g).
\end{eqnarray*}
Let $D(\underline{v})$ be the image of the map $j_{\underline{v}}$.
Then $D(\underline{v})\simeq
I(\underline{v})/I(\underline{v}+\underline{1})$ and we define
$F_{\underline{v}}$ as $D(\underline{v}) \cap
(D_1^{*}(\underline{v}) \times \cdots \times
D_r^{*}(\underline{v}))$ where $D_j^{*}(\underline{v})$ denotes
$D_j(\underline{v}) \setminus \{\underline{0}\}$, $j \in J$.
Having the map
\begin{eqnarray*}
\rho: \quad\hat{S}_{H} \quad & \longrightarrow & S_{H} \\
(\underline{\nu}(g),\underline{a}(g)) & \longmapsto &
\underline{\nu}(g),
\end{eqnarray*}
$F_{\underline{v}}$ can also be expressed as
$\rho^{-1}(\underline{v})$ and therefore one also calls the spaces
$F_{\underline{v}}$ the fibers of the extended semigroup
$\hat{S}_{H}$. For $\underline{v} \in S$, the space
$F_{\underline{v}}$ is the complement to an arrangement of vector
subspaces in the vector space $D(\underline{v})$. Moreover
$F_{\underline{v}}$ is invariant with respect to multiplication by
nonzero constants. Let
$\mathbb{P}F_{\underline{v}}=F_{\underline{v}}/\mathbb{C}^{*}$ be
the projectivization of $F_{\underline{v}}$. It is the complement to
an arrangement of projective subspaces in the projective space
$\mathbb{P}D(\underline{v})$. \\As the ideal $I(\underline{v})$ is a
monomial ideal, one can see the space $F_{\underline{v}}$ as the set
of functions $g=\sum_{i=1}^{s}\lambda_i m_i$ in $\mathcal{O}_{X,o}$
(the $m_i$ are monomials and the $\lambda_i, 1 \leq i \leq s$, are
complex numbers different from $0$) for which
$\underline{\nu}(g)=\underline{v}$ and for which for all $i \in
\{1,\cdots,s\}$ holds that $m_i \in I(\underline{v})$ and that there
exists a $j \in J$ such that $\nu_j(m_i)=\nu_j(g)$. If a function
$g$ with $\underline{\nu}(g)=\underline{v}$ has this form, we say
that $g$ is in reduced form. We write supp($g$) for the support of
$g$, which is the set $\{m_i \mid 1 \leq i \leq s\}$.
\\ Now let us fix $\underline{v} \in \mathbb{Z}^r$. Let $\mathcal{M}$ be the set of
all monomials that can appear in the support of some $g$ in
$\mathbb{P}F_{\underline{v}}$, so $\mathcal{M}=\{m \mbox{ monomial
$|$ } m \in I(\underline{v}) \mbox{ and } \exists j \in J:
\nu_j(m)=v_j\}$. Note that $\mathcal{M}$ is a finite set. For a
subset $L$ of $\mathcal{M}$, let $\underline{\nu}(L)$ be the vector
$\underline{w} \in \mathbb{Z}^r$ with $w_j=\mbox{min}\{\nu_j(m) \mid
m \in L\}$, $j \in J$.
\newpage \noindent \textbf{Proposition 1.---}
\begin{center} $P(\underline{t})=\sum_{\underline{v} \in
\mathbb{Z}^r}\chi(\mathbb{P}F_{\underline{v}})\underline{t}^{\underline{v}}$.
\end{center}
\emph{Proof.}
\\We write $\mathbb{P}F_{\underline{v}}$ as a disjoint union:
\[\mathbb{P}F_{\underline{v}}=\bigcup_{L \subset \mathcal{M}, \underline{\nu}(L)=\underline{v}}\{g \in \mathbb{P}F_{\underline{v}} \mbox{ $|$ supp($g$)}=L\}
.\] For $\Lambda_L=\{g \in \mathbb{P}F_{\underline{v}} \mbox{ $|$
supp($g$)}=L\}$, one has that
$\chi(\Lambda_L)=\chi({\mathbb{C}^{*}}^k)=0$ for some $k \in
\mathbb{Z}_{>0}$ when $L$ is not a singleton, and
$\chi(\Lambda_L)=1$ when $L$ is a singleton. This gives us
\[\chi(\mathbb{P}F_{\underline{v}})=N(\underline{v}).\] \hfill $\blacksquare$
\\ \\ This result can also be obtained by the combinatorial proof
given in \cite{trio2} or in \cite{trio3}. Actually their proof works
for every variety for which the Poincar\'{e} series exists.
Exploiting the fact that our varieties are toric, we obtain this
explicit and much shorter proof.
\\ \\ \indent Analogously to the case of curves and rational surface singularities,
we write the Poincar\'{e} series as an integral with respect to the
Euler characteristic over the projectivization
$\mathbb{P}\mathcal{O}_{X,o}$. We first recall this notion,
introduced in \cite{trio3} and inspired by the notion of motivic
integration (see for example \cite{denefloeser}). It was developed
to integrate over $\mathbb{P}\mathcal{O}_{\mathbb{C}^n,o}$ what is
not allowed by the usual Viro construction where one integrates with
respect to the Euler characteristic over finite dimensional spaces
(see \cite{Viro}). It can be extended to integrals over
$\mathbb{P}\mathcal{O}_{X,o}$, for $X$ an arbitrary variety.
\\ \indent Let $\mathsf{m}$ be the maximal ideal in the
ring $\mathcal{O}_{X,o}$ and for $k \in \mathbb{Z}_{\geq 0}$, let
$\mathcal{J}^{k}_{X,o}=\mathcal{O}_{X,o}/\mathsf{m}^{k+1}$ be the
space of $k$-jets of functions on the toric variety $X$. For a
complex vector space $L$ (finite of infinite dimensional) let
$\mathbb{P}L=(L \setminus \{0\})/\mathbb{C}^*$ be its
projectivization and let $\mathbb{P}^*L$ be the disjoint union of
$\mathbb{P}L$ with a point. One has a natural map
$\pi_k: \mathbb{P}\mathcal{O}_{X,o} \mapsto
\mathbb{P}^*\mathcal{J}^{k}_{X,o}$.
\\ \\
\noindent \textbf{Definition.---} A subset $A \subset
\mathbb{P}\mathcal{O}_{X,o}$ is said to be cylindrical if
$A=\pi_k^{-1}(B)$ for a constructible (i.e. a finite union of
locally closed sets) subset $B \subset
\mathbb{P}\mathcal{J}^{k}_{X,o} \subset
\mathbb{P}^*\mathcal{J}^{k}_{X,o}$. \\ \\
\textbf{Definition.---} For a cylindrical subset $A \subset
\mathbb{P}\mathcal{O}_{X,o}$ ($A=\pi_k^{-1}(B)$, $B \subset
\mathbb{P}\mathcal{J}^{k}_{X,o}$), its Euler characteristic
$\chi(A)$ is defined as the Euler characteristic $\chi(B)$ of the
set $B$. \\ \\ Let $\psi: \mathbb{P}\mathcal{O}_{X,o} \rightarrow G$
be a function which takes values in an abelian group $G$.
\\ \\
\textbf{Definition.---} We say that the function $\psi$ is
cylindrical if, for each $g \in G, g \neq 0$, the set $\psi^{-1}(g)
\subset \mathbb{P}\mathcal{O}_{X,o}$ is cylindrical. \newpage
\noindent \textbf{Definition.---} The integral of a cylindrical
function $\psi$ over the space $\mathbb{P}\mathcal{O}_{X,o}$ with
respect to the Euler characteristic is
\[\int_{\mathbb{P}\mathcal{O}_{X,o}} \psi d\chi := \sum_{g \in G, g \neq 0} \chi(\psi^{-1}(g))\cdot g\]
if this sum makes sense in $G$. If the integral exists, the function
$\psi$ is said to be integrable.\\ \indent 
Let $\mathbb{Z}[[\underline{t}]]=\mathbb{Z}[[t_1,\cdots,t_r]]$ be
the group with respect to the addition of formal power series in the
variables $t_1,\cdots,t_r$. We have the map
$\underline{\nu}: \mathbb{P}\mathcal{O}_{X,o} \longrightarrow
\mathbb{Z}^r$. 
Let $\underline{t}^{\underline{\nu}}$ be the corresponding function
with values in $\mathbb{Z}[[\underline{t}]]$. 
\\ \\ \textbf{Proposition 2.---}
\begin{center}
$P(\underline{t})=\int_{\mathbb{P}\mathcal{O}_{X,o}}\underline{t}^{\underline{\nu}}d\chi.$
\end{center}
\emph{Proof.}
For $\underline{v} \in \mathbb{Z}_{\geq 0}^r$, let $N=1 +$ max$\{v_i
\mid 1 \leq i \leq r\}$ and let $Y_{\underline{v}}= \{j^N g \in
\mathbb{P}\mathcal{J}^{N}_{X,o} \mid g \in
\mathbb{P}F_{\underline{v}}\} \subset
\mathbb{P}\mathcal{J}^{N}_{X,o}$. Then $\{g \in \mathbb{P}
\mathcal{O}_{X,o} \mid
\underline{\nu}(g)=\underline{v}\}=\pi_N^{-1}(Y_{\underline{v}})$.
Consider the map
\begin{eqnarray*}
\alpha: Y_{\underline{v}} & \longrightarrow &
F_{\underline{v}}\ \\
j^N g & \longmapsto & \underline{a}(g).
\end{eqnarray*}
This map is $\mathbb{C}^*$-invariant and so can be considered as a
map $\alpha: Y_{\underline{v}}  \longrightarrow
\mathbb{P}F_{\underline{v}}$. As $\alpha$ is a locally trivial
fibration whose fibre is a complex affine space we obtain
\[\int_{\mathbb{P}\mathcal{O}_{X,o}}\underline{t}^{\underline{\nu}}d\chi=
\sum_{\underline{v} \in \mathbb{Z}^r}
\chi(Y_{\underline{v}})\underline{t}^{\underline{v}}=
\sum_{\underline{v} \in \mathbb{Z}^r}
\chi(\mathbb{P}F_{\underline{v}})\underline{t}^{\underline{v}}=P(\underline{t}).
\]\hfill
$\blacksquare$
 \begin{center} \textsc{2. Poincar\'{e} series
of $\mathbb{C}^d$ induced by a toric constellation}
\end{center} ${}$\\
\indent In this section we look at the particular case where $X$ is
$\mathbb{C}^d$ endowed with the action of the torus $T\cong
(\mathbb{C}^*)^d$ ($d \geq 2$). We study the Poincar\'{e} series of
$\mathbb{C}^d$ where the modification $\pi$ is given by a toric
constellation $\mathcal{C}$ with origin. Let us first recall the
basic notions in the context of toric constellations (see for
example \cite{clusters}).\\ \indent Let $Z$ be a variety obtained
from $X$ by a finite succession of point blowing-ups. A point $Q \in
Z$ is said to be infinitely near to a point $O \in X$ if $O$ is in
the image of $Q$; we write $Q \geq O$. A constellation is a finite
sequence $\mathcal{C}=\{Q_0,Q_1,\cdots,Q_{r-1}\}$ of infinitely near
points of $X$ such that $Q_0 \in X=X_0$ and each $Q_j$ is a point on
the variety $X_j$ obtained by blowing up $Q_{j-1}$ in $X_{j-1}$, $j
\in J$. When $Q_j$ is a $0$-dimensional $T$-orbit in the toric
variety $X_j$, for $j \in J$, the constellation is said to be toric. \\
The relation `$\geq$' gives rise to a partial ordering on the points
of a constellation. In the case that they are totally ordered, so
$Q_r \geq \cdots \geq Q_0$, the constellation $\mathcal{C}$ is
called a chain.
\\For every $Q_j$ in $\mathcal{C}$, the subsequence $C^j=\{Q_i \mid Q_j \geq
Q_i\}$ of $\mathcal{C}$ is a chain. The integer $l(Q_j)=\#
\mathcal{C}^j - 1$ is called the level of $Q_j$. In particular $Q_0$
has level $0$. If no other point of $\mathcal{C}$ has level $0$ then
$Q_0$ is called the origin of $\mathcal{C}$. \\A tree with a root
such that each vertex has at most $d$ following adjacent vertices is
called a $d$-ary tree. There is a natural bijection between the set
of $d$-dimensional toric constellations with origin and the set of
finite $d$-ary trees with a root, with the edges weighted with
positive integers not greater than $d$, such that two edges with the
same source have different weights. We simply say `weights' when we
want to refer to the weights on the edges. The weights indicate in
which affine chart the points of the constellation are created (see
\cite{clusters} for details).
\\ \\
\indent Let $E_j$ ($j \in J$) be the irreducible components of the
exceptional divisor $\mathcal{D}$ created by blowing up the
constellation $\mathcal{C}=\{Q_0,Q_1,\cdots,Q_{r-1}\}$ and let
$\xi_j$ be the generic point of $E_j$ ($j \in J$). Then
$\mathcal{O}_{X,\xi_j}$ is a discrete valuation ring. We denote the
induced valuation by $\nu_j$. \\We denote the matrix of the linear
system of equations $\langle s,\underline{\nu}\rangle=\underline{v}$
by $\mathcal{L}(\mathcal{C})$ and we denote the column vectors of
$\mathcal{L}(\mathcal{C})$ by $v_1, \cdots, v_d$. Let $C$ be the
cone in $\mathbb{Z}^r_{\geq 0}$ generated by $v_1, \cdots, v_d$.
Note that $C$ is the monomial cone associated to the constellation
$\mathcal{C}$. Recall that the cone is regular if it can be
generated by a part of a basis of $\mathbb{Z}^r$. If
$\mathcal{L}(\mathcal{C})$ has rank $s$ and if $s < d$, then the
cone is said to be degenerated.
\\ \\ \textbf{Example:} \\
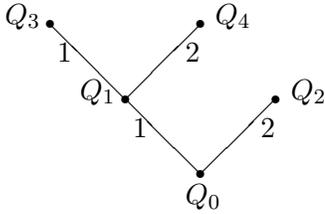
\begin{tabular}{p{4.5cm}p{8cm}}
\begin{picture}(35,10)(-8,-2)
\put(20,-25){\circle*{1}}
\put(30,-15){\circle*{1}}\put(10,-15){\circle*{1}}
\put(20,-25){\line(1,1){10}} \put(20,-25){\line(-1,1){10}}
\put(18,-29){$Q_0$} \put(32,-15){$Q_2$}  \put(4,-15){$Q_1$}
\put(28,-20){$2$} \put(11,-20){$1$} \put(10,-15){\line(1,1){10}}
\put(10,-15){\line(-1,1){10}}\put(20,-5){\circle*{1}}\put(0,-5){\circle*{1}}
\put(22,-5){$Q_4$}  \put(-6,-5){$Q_3$} \put(18,-10){$2$}
\put(1,-10){$1$}
\end{picture}
& Suppose $d=4$ and $\mathcal{C}$ is the constellation pictured at
the left. By blowing up in the origin $Q_0$ we get an exceptional
variety $B_0 \cong \mathbb{P}^3$. In $B_0$ there are $2$ points in
which we blow up, namely $Q_1$ and $Q_2$. For example the point
$Q_1$ is the origin of the affine chart induced by the weight $1$,
we shortly denote this affine chart by `$1$'.
\end{tabular}
\\ After blowing
up in $Q_1$ we get an exceptional variety $B_1 \cong \mathbb{P}^3$,
where again we blow up in $2$ points. The point $Q_3$ is the origin
of the affine chart `$1-1$' and $Q_4$ is the origin of the affine
chart `$1-2$'. The associated linear system
$\mathcal{L}(\mathcal{C})$ is \\ \\$\left\{
\begin{array}{ccccccccc} a & + & b & + & c & + &
d & = & v_1 \\
a & + & 2b & + & 2c & + &
2d & = & v_2 \\
2a & + & b & + & 2c & + &
2d & = & v_3 \\
a & + & 3b & + & 3c & + &
3d & = & v_4 \\
2a & + & 3b & + & 4c & + & 4d & = & v_5
\end{array} \right.$
\\ \\
and $C$ is the cone $\langle
(1,1,2,1,2),(1,2,1,3,3),(1,2,2,3,4)\rangle \subset
\mathbb{Z}^5_{\geq 0}$ which is obviously degenerated.
\vspace*{-0.5cm} \[\]\hfill $\square$ \newpage \indent To know the
Poincar\'{e} series one can use Theorem $1$. Let $\underline{v}_1,
\cdots, \underline{v}_{d}$ be the column vectors of
$\mathcal{L}(\mathcal{C})$. Then it follows by Theorem $1$ that
\[P(\underline{t})=\frac{1}{(1-\underline{t}^{\underline{v}_1})\cdots (1-\underline{t}^{\underline{v}_d})}.\]
One can also obtain $P(\underline{t})$ by computing the numbers
$N(\underline{v})=\#\{s \in \mathbb{N}^d \mid \langle
s,\underline{\nu}\rangle=\underline{v}\}$ for each $\underline{v}
\in \mathbb{Z}^r$. We determine them as later these values will be
useful for us. The following proposition gives some properties of
the cone $C$ which will allow us to compute the numbers
$N(\underline{v})$.
\\ \\
\noindent \textbf{Proposition 3.---}
\begin{enumerate}
\item \emph{$C$ is a regular cone;} \item \emph{$C$ is degenerated
if and only if the number of different weights appearing in the
constellation is less than or equal to $d-2$.}
\end{enumerate}
\emph{Proof.} If the number of different weights appearing in the
constellation $\mathcal{C}$ is less than or equal to $d-2$, then at
least $2$ columns in $\mathcal{L}(\mathcal{C})$ are equal and $C$ is
degenerated. Suppose that the number of different weights is bigger
than $d-2$, say that $1,\cdots,d-1$ are weights appearing in the
constellation and let $Q_1,\cdots,Q_{d-1}$ be points in the
constellation such that $Q_i$ is a point with minimal level arising
in an affine chart induced by the weight $i$ and such that
whenever $Q_i \geq Q_j$ and $Q_i \neq Q_j$, then $i > j$ ($1 \leq i, j \leq d-1$). \\
The linear equations induced by the origin $Q_0$ and by $Q_1,
\cdots, Q_{d-1}$ give rise to a linear system whose determinant can
be supposed to be of the form
\[\begin{tabular}{|c c c c c c|}
$1$ & $1$ & $1$ & $\cdots$ & $\cdots$ & $1$ \\
$1$ & $2$ & $2$ & $\cdots$ & $\cdots$ & $2$ \\
$*$ & $a_2-1$ & $a_2$ & $\cdots$ & $\cdots$ & $a_2$\\
$*$ & $*$ & $a_3-1$ & $a_3$ & $\cdots$ & $a_3$\\
$\vdots$ & $\vdots$ & $\vdots$ & $\vdots$ & $\vdots$ & $\vdots$ \\
$*$ & $\cdots$ & $\cdots$ & $*$ & $a_{d-1}-1$ & $a_{d-1}$
\end{tabular}=
\begin{tabular}{|c c c c c c|}
$1$ & $1$ & $1$ & $\cdots$ & $\cdots$ & $1$ \\
$-1$ & $0$ & $0$ & $\cdots$ & $\cdots$ & $0$ \\
$*$ & $-1$ & $0$ & $\cdots$ & $\cdots$ & $0$\\
$*$ & $*$ & $-1$ & $0$ & $\cdots$ & $0$\\
$\vdots$ & $\vdots$ & $\vdots$ & $\vdots$ & $\vdots$ & $\vdots$ \\
$*$ & $\cdots$ & $\cdots$ & $*$ & $-1$ & $0$
\end{tabular},\]
where $a_2, \cdots, a_{d-1}$ are integer numbers. As this
determinant is equal to $1$, the cone $C$ is regular and non-degenerated.\\
Now if $C$ is degenerated, and generated by
$s=\mbox{rank}(\mathcal{L}(\mathcal{C}))$ vectors then there are
$s-1$ different weights appearing in the constellation. In the same
way as above we obtain a $(s \times s$)-determinant which is equal
to $1$ such that also in this case $C$ is regular.
\vspace*{-0.5cm}\[\] \hfill $\blacksquare$ \\ \\
\noindent Note that it follows from the proof that saying that the
cone $C$ is degenerated, is the same as saying that there is one
column that appears at least twice in $\mathcal{L}(\mathcal{C})$ and
that there are no other linear dependencies between the columns of
$\mathcal{L}(\mathcal{C})$.
\begin{itemize}
\item If $C$ is non-degenerated then $N(\underline{v})=1$ if $\underline{v} \in
C$, else $N(\underline{v})=0$. Suppose that $\underline{v}_1,
\cdots, \underline{v}_d$ are the column vectors of the linear system
$\mathcal{L}(\mathcal{C})$. As $C$ is regular, we obtain again
\[P(\underline{t})=\frac{1}{(1-\underline{t}^{\underline{v}_1})\cdots (1-\underline{t}^{\underline{v}_d})}.\]
\item When $C$ is degenerated, let then $\underline{v}_1, \cdots,
\underline{v}_{s}$ be the $s$ different vectors that generate $C$
and let $\underline{v}_{s}$ be the vector that appears at least
twice as column vector in $\mathcal{L}(\mathcal{C})$. As $C$ is
regular, one can write each $\underline{v} \in C$ in a unique way
as $\underline{v}=\lambda_1 \underline{v}_1 + \cdots +
\lambda_{s-1} \underline{v}_{s-1} + \lambda \underline{v}_{s}$,
for some $\lambda_i, \lambda \in \mathbb{Z}_{\geq 0}$ ($1 \leq i
\leq s-1$).\\ Setting $k=d-s+1$, a simple calculation shows that
\[N(\underline{v})=\begin{pmatrix}k + \lambda - 1\\\lambda
\end{pmatrix}\] if $\underline{v} \in
C$, else $N(\underline{v})=0$. Using that $\sum_{\lambda \in
\mathbb{Z}_{\geq 0}}\begin{pmatrix}k + \lambda - 1 \\ \lambda
\end{pmatrix}x^{\lambda}= \frac{1}{(1-x)^k}$, one sees again that
\[P(\underline{t})=\frac{1}{(1-\underline{t}^{\underline{v}_1})\cdots (1-\underline{t}^{\underline{v}_{s-1}})(1-\underline{t}^{\underline{v}_{s}})^k}.\]

\end{itemize}
\noindent \\ \\ For the example given above the Poincar\'{e} series
is then
\[\frac{1}{(1-t_1t_2t_3^2t_4t_5^2)
(1-t_1t_2^2t_3t_4^3t_5^3)(1-t_1t_2^2t_3^2t_4^3t_5^4)^2}.\] \\ \\
\indent For curves, rational surface singularities and plane
divisorial valuations, there exists a description of the
Poincar\'{e} series at the level of the modification space. Let
$\mathcal{D}=\bigcup_{j=1}^r E_j$ be the exceptional variety with
irreducible components $E_j, j \in J$. We denote by
$\overset{\circ}{E_j}$ the smooth part of the irreducible component
$E_j$, i.e. without intersection points with all other components of
the exceptional divisor. Let $M=-(E_i \circ E_j)$ be minus the
intersection matrix of the components of the exceptional variety
$\mathcal{D}$. Let $\nu_j$ be the discrete valuation on the local
ring $\mathcal{O}_{X,o}$ induced by $E_j$. The semigroup of values
$S:=\{\underline{\nu}(g) \mid g \in \mathcal{O}_{X,o}\}$ is exactly
the set of vectors $\{\underline{v} \in \mathbb{Z}^r_{\geq 0} \mid
\underline{v}M \geq \underline{0}\}$. For a topological space $E$,
let $S^nE=E^n/S_n$ ($n \geq 0$) be the $n$-th symmetric power of the
space $E$, i.e. the space of $n$-tuples of points of the space $E$
($S^0E$ is a point).
\\Campillo, Delgado and Gusein-Zade construct the space
\[Y=\bigcup_{\{\underline{v} \in S\}}\left(\prod_{j=1}^r S^{n_j}\overset{\circ}{E_j}\right),\]
where
$\underline{v}M=(n_1,\cdots,n_r)=:\underline{n}(\underline{v})$. For
$g \in \mathcal{O}_{X,o}$, $g\neq 0$ and
$\underline{v}=\underline{\nu}(g)$, the number $n_j(\underline{v})$
is equal to the intersection number of the strict transform of $g$
with $E_j$. Let $Y_{\underline{v}}$ be the connected component
$\prod_{j=1}^r S^{n_j}\overset{\circ}{E_j}$ of $Y$, where
$(n_1,\cdots,n_r)=\underline{n}(\underline{v})$. They show that
\begin{eqnarray}
P(\underline{t})=\sum_{\underline{v} \in \mathbb{Z}^r}
\chi(Y_{\underline{v}})\underline{t}^{\underline{v}}.
\end{eqnarray}
In the case of curves and plane divisorial valuations, this
description induces the elegant formula of the Poincar\'{e} series
where the exponents can be written as the Euler characteristics of
the smooth parts $\overset{\circ}{E_j}$.
\\ \\ \indent In what follows we prove a generalized form of $(1)$
for toric constellations. I want to thank Sabir Gusein-Zade for the
very stimulating conversation about this subject. \\For $j \in J$,
let $B_j$ be the projective $(d-1)$-dimensional exceptional variety
that is created by blowing up the point $Q_j$ in $X_j$. When we want
to refer to the strict transform of $B_j$ at some intermediate
stadium, we also write $E_j$. If $g \in
\mathcal{O}_{\mathbb{C}^d,o}$, then we write $\hat{g}$ for the
strict transform at the end of the process as well
as for the strict transforms of $g$ in the intermediate stadia. 
For $v \in \mathbb{Z}^r$, we define the set
\begin{eqnarray*}
D_{\underline{v}}: & = & \{ \{\hat{g}=0\} \cap \mathcal{D} \mid g
\in \mathcal{O}_{\mathbb{C}^d,o}, \mbox{
$\underline{\nu}(g)=\underline{v}$ and } \{\hat{g}=0\} \mbox{ does
not contain}
\\ & & \mbox{any non-empty intersection } E_a \cap E_b, a, b \in J, a \neq
b\}.
\end{eqnarray*}
Obviously to know $D_{\underline{v}}$ it is sufficient to consider
the elements $g$ in $\mathbb{P}F_{\underline{v}}$. We make a
topological space of it as follows. We write $E_{\underline{v}}$ for
the sum $\sum_{i=1}^r v_iE_i$ and $M$ for the line bundle associated
to $E_{\underline{v}}$. The restriction $R$ of
$\mathcal{O}_Z(-E_{\underline{v}}) \otimes M^{-1}$ to $\mathcal{D}$
is a line bundle and as $\mathcal{D}$ is a projective variety, the
global sections of $R$ form a finite dimensional vector space. For
$g \in F_{\underline{v}}$, the divisor $\hat{g} \cap \mathcal{D}$ is
the divisor of zeroes of a global section of $R$. Then
$D_{\underline{v}}$ can be seen as a subset of the projectivization
of this vector space.
\\ \\
\textbf{Theorem 2.---} \emph{The Poincar\'{e} series
P(\underline{t}) is equal to}
\[\sum_{\underline{v} \in \mathbb{Z}^r}
\chi(D_{\underline{v}})\underline{t}^{\underline{v}}.\] To achieve
the result we will construct a subspace $Z_{\underline{v}}$ of
$\mathbb{P}F_{\underline{v}}$ that has the same Euler characteristic
as $\mathbb{P}F_{\underline{v}}$ and such that there exists a
homeomorphism of $Z_{\underline{v}}$ with $D_{\underline{v}}$. Then Theorem $2$ will follow by Proposition $1$.\\
The obvious candidate for the space $Z_{\underline{v}}$ is the set
\begin{eqnarray*}  Z_{\underline{v}} & := & \{g \in \mathbb{P}F_{\underline{v}} \mid \{\hat{g}=0\} \mbox{ does not
contain any non-empty} \\
& & \mbox{intersection } E_a \cap E_b, a, b \in J, a \neq b\}.
\end{eqnarray*}\newpage
\noindent \textbf{Lemma 1.---}
\[\chi(Z_{\underline{v}})=\chi(\mathbb{P}F_{\underline{v}}).\]
\emph{Proof.} Let $g=\sum_{i=1}^{s}\lambda_i m_i$ be in reduced form
and suppose that $E_a \cap E_b \neq \emptyset$ for some $a$ and $b$
in $J$, $a \neq b$. If $g$ contains $E_a \cap E_b$ then also
$g_{\underline{\mu}}=\sum_{i=1}^{s}\mu_i m_i$ contains $E_a \cap
E_b$, for all $\mu=(\mu_1,\cdots,\mu_s) \in {\mathbb{C}^*}^s$. This
yields that \begin{eqnarray*} \chi(Z_{\underline{v}}) & = & \chi(\{g
\in \mathbb{P}F_{\underline{v}} \mid g \mbox{ a
monomial and } \{\hat{g}=0\} \mbox{ does not contain any} \\
& &
\mbox{non-empty intersection } E_a \cap E_b,  a, b \in J, a \neq b\})\\
& = & \chi(\{g \in \mathbb{P}F_{\underline{v}} \mid g \mbox{ a
monomial}\}) \\
& = & N(\underline{v}) \\
& = & \chi(\mathbb{P}F_{\underline{v}}).
\end{eqnarray*}
 \hfill $\blacksquare$
\\ \\
\noindent Now we investigate the map
\begin{eqnarray*}
\phi: Z_{\underline{v}} & \longrightarrow & D_{\underline{v}}\\
 g & \longmapsto & \{\hat{g}=0\} \cap \mathcal{D}.
\end{eqnarray*}
The following lemma tells us how $\{\hat{g}=0\} \cap \mathcal{D}$ looks like.\\ \\
\textbf{Lemma 2.---} \emph{Consider $g \in
\mathbb{P}F_{\underline{v}}$, a function in reduced form, say
$g=\sum_{i=1}^{s}\lambda_i m_i$. For $j \in J$, let $\Lambda_{g,j}$
be the set $\{m \in \mbox{ supp}(g) \mid \nu_j(m)=v_j\}$. Then the
equation of $\{\hat{g}=0\} \cap E_j$ is $\sum_{m_i \in
\Lambda_{g,j}}\lambda_i \hat{m}_i=0$ in every affine chart where the
intersection $\{\hat{g}=0\} \cap E_j$ is visible.}
\\ \\\emph{Proof.} When $\{\hat{g}=0\} \cap B_j \neq \emptyset$, then in an affine
chart covering $B_j \cong \mathbb{P}^{d-1}$ it can be described by
the equation $\sum_{i \in K_j}\lambda_i \hat{m}_i=0$ where
\begin{eqnarray*}
K_j & = & \{i \mid \nu_j(m_i)=\mbox{min}\{\nu_j(m) \mid m \in
\mbox{supp}(g)\}\}\\
&= & \{i \mid \nu_j(m_i)=v_j\} \\
& = & \Lambda_{g,j}.
\end{eqnarray*}
Now Lemma 2 follows directly. \vspace*{-0.5cm}\[\]\hfill
$\blacksquare$
\\ \\
Note that $\hat{m}$ can be equal to $1$ in some affine chart
covering $E_j$ ($j \in J$). However, there exists always a $j \in J$
such that $\{\hat{m}=0\} \cap E_j \neq \emptyset$.
\\ \\ Let us have a look at the following example. Lemma $2$ allows us to deduce quickly
a necessary condition on a subset $S$ of
$\mathbb{P}F_{\underline{v}}$ for the sets $\{g \in S\}$ and
$\{\{\hat{g}=0\} \cap \mathcal{D} \mid g \in S \}$ to be in $1-1$
correspondence. \newpage \noindent \textbf{Example (continued).---} \\ \\
Let $\underline{v}$ be $(6,11,10,14,18)$. Take
$g(x,y,z,u)=y^2z^4+x^3y^4+x^{14}$. The values of the monomials in
the support of $g$ are:
\[\nu(y^2z^4)=(6,12,10,18,22), \quad \nu(x^3y^4)=(7,11,10,15,18) \quad \mbox{ and}\]
\[\nu(x^{14})=(14,14,28,14,28).\] Lemma $2$ tells us that the
equation of $\hat{g} \cap E_1$ in an affine chart where the
intersection is visible, is $\hat{y^2z^4}$. For $\hat{g} \cap E_2$
it is $\hat{x^3y^4}$, for $\hat{g} \cap E_3$ it is
$\hat{y^2z^4}+\hat{x^3y^4}$, for $\hat{g} \cap E_4$ it is
$\hat{x^{14}}$ and for $\hat{g} \cap E_5$ it is $\hat{x^3y^4}$. It
follows that the strict transform of $h(x,y,z,u)=y^2z^4+x^3y^4+ \mu
x^{14}$, with $\mu \neq 0$, has the same intersection with
$\mathcal{D}$ as $\hat{g}$. \vspace*{-1cm}\[\] \hfill $\square$
\\ \\
\noindent We formalize what the example shows. As in Proposition
$1$, write $\mathbb{P}F_{\underline{v}}$ as the disjoint union
$\bigcup_{L \subset \mathcal{M}, \underline{\nu}(L)=\underline{v}}
\Lambda_L$. Let $L$ be a support appearing in this disjoint union.
For $j \in J$, we define $\Lambda_{L,j}$ as the set of monomials $m$
in $L$ such that $\nu_j(m)=v_j$. The observation made in the example
shows that the map
\begin{eqnarray*}
\phi_L: \{g \in \mathbb{P}F_{\underline{v}}\mid \mbox{supp}(g)=L\} &
\longrightarrow & \{\{\hat{g}=0\} \cap \mathcal{D} \mid g \in
\mathbb{P}F_{\underline{v}} \mbox{ and supp}(g)=L\} \\
g & \longmapsto & \{\hat{g}=0\} \cap \mathcal{D}
\end{eqnarray*}
certainly can not be a bijection if there exists a subset $D$ of
$J$, with the following properties:
\begin{eqnarray}
\exists a, b \in J: a \in D, \quad b \notin D, \quad \left(\cup_{d
\in D}\Lambda_{L,d}\right) \cap \left(\cup_{d \notin
D}\Lambda_{L,d}\right)=\emptyset.
\end{eqnarray}Indeed, if such a subset $D$ exists, then write
$g=g_a+g_b$ where supp$(g_a)=\cup_{d \in D}\Lambda_{L,d}$ and
supp$(g_b)=\cup_{d \notin D}\Lambda_{L,d}$. Then by Lemma $2$ it
follows that for
$g_{\underline{\lambda}}=\lambda_ag_a+\lambda_bg_b$, the transforms
$\{\hat{g_{\underline{\lambda}}}=0\}$ and $\{\hat{g}=0\}$ have the
same intersection with $\mathcal{D}$ for all
$\underline{\lambda}=(\lambda_a,\lambda_b) \in {\mathbb{C}^*}^2$. We
claim that also the converse is true.
\\ \\ \noindent
\textbf{Proposition 4.---} \emph{Let $L \subset \mathcal{M}$ be as
above. If for all $a, b \in J$, $a \neq b$, holds that there exists
no subset $D$ with $a \in D$, $b \notin D$ and $\left(\cup_{d \in
D}\Lambda_{L,d}\right) \cap \left(\cup_{d \notin
D}\Lambda_{L,d}\right)=\emptyset$, then the map $\phi_L$ is a
bijection.}
\\ \\ \emph{Proof.}
\begin{enumerate}
\item First we show that for a monomial $m$ in $L$, when given $\{\hat{m}=0\} \cap \mathcal{D}$, one can find $m$ again.
So suppose $m$ in $L$. Take a $j \in J$ such that $\nu_j(m)=v_j$ and
such that $\{\hat{m}=0\} \cap E_j$ is visible in an affine chart in
the final stadium, say in the chart presented as $c_1- c_2-\cdots-
c_t$. Take the subchain of the constellation with consecutive
weights $c_1, c_2, \cdots, c_{t-1}$ and add an edge with weight
$c_t$ at the top. We call this new chain $K$. For $i \in
\{1,\cdots,d\}$, let $Q_{K_i}$ be the point with maximal level of
$K$ for which the weight going out from $Q_{K_i}$ is $i$, if this
point exists. If $m=x_1^{n_1}x_2^{n_2}\cdots x_d^{n_d}$, then the
equation of $\hat{m}$ in the considered chart is
\[x_1^{f_1(n_1,\cdots,n_d)-h_1}x_2^{f_2(n_1,\cdots,n_d)-h_2}\cdots x_d^{f_d(n_1,\cdots,n_d)-h_d},\]
where $f_i$ is the left hand side of the linear equation induced by
the point $Q_{K_i}$ in $\mathcal{L}(\mathcal{C})$ if $Q_{K_i}$
exists and $f_i(n_1,\cdots,n_d)=n_i$ if this point is not defined.
The value $h_i$ is equal to $v_{K_i}$ if $Q_{K_i}$ exists and else
$h_i=0$.
\\Now suppose that $\hat{m} \cap E_i=x_1^{n'_1}x_2^{n'_2}\cdots x_d^{n'_d}$ is given. When $Q_j$
exists for every $j \in \{1,\cdots,d\}$, it follows from Proposition
$3$ that the linear system $\{f_j(\underline{n})-h_j=n'_j\}_{1 \leq
j \leq d}$ has a \emph{unique} solution for $\underline{n}$. If not
all points $Q_j$ are defined, one also sees immediately that the
linear system $\{f_j(\underline{n})-h_j=n'_j\}_{1 \leq j \leq d}$ is
regular. \item Let $g$ and $h$ be two functions in
$\mathbb{P}F_{\underline{v}}$ with the same support $L$. Write
$g=\sum_{i=1}^{s}\lambda_i m_i$ and $h=\sum_{i=1}^{s}\mu_i m_i$,
with $\lambda_i$ and $\mu_i$ different from $0$ ($1 \leq i \leq s$).
Suppose that $\{\hat{g}=0\} \cap \mathcal{D}=\{\hat{h}=0\} \cap
\mathcal{D}$ and that $\lambda_j \neq \mu_j$. For lack of a subset
$D$ of $J$ with property $(2)$ and because of part $1$ of the proof,
it follows that $\lambda_i/\mu_i=\lambda_j/\mu_j$, for all $i \in J$
and so $g=h$.
\end{enumerate}
\hfill $\blacksquare$
\\ \\
\indent The functions we are interested in are the functions $g$
such that $\{\hat{g}=0\}$ does not contain any non-empty
intersection $E_a \cap E_b$.
They can be characterized as follows:\\ \\
\textbf{Lemma 3.---}
\emph{Let $a$ and $b$ be different elements of $J$ such that $E_a
\cap E_b \neq \emptyset$. Then
\begin{center}
$\{\hat{g}=0\}$ contains $E_a \cap E_b$ \\
$\Updownarrow$\\
there is no $m$ in $\mbox{supp}(g)$ for which $\nu_a(m)=v_a$ and
$\nu_b(m)=v_b$.
\end{center}}
\noindent \emph{Proof.}
\\ \indent Suppose that there is no $m$ in $\mbox{supp}(g)$ for which $\nu_a(m)=v_a$ and
$\nu_b(m)=v_b$. In an affine chart where one sees $E_a \cap E_b$,
one has \[\hat{g}=\hat{g_a}+\hat{g}_b+\hat{r}, \quad E_a
\leftrightarrow x_a=0, \quad \quad E_b \leftrightarrow x_b=0\] where
$g_a$ is the part of $g$ with supp$(g_a)=\{m \in \mbox{supp}(g) \mid
\nu_a(m)=\nu_a(g)\}$, $g_b$ is the part of $g$ with supp$(g_b)=\{m
\in \mbox{supp}(g) \mid \nu_b(m)=\nu_b(g)\}$ and $r$ is $g-g_a-g_b$.
From Lemma $2$ it follows that $\hat{g_b} + \hat{r} \in (x_a)$ and
$\hat{g_a} + \hat{r} \in (x_b)$. Then also $\hat{g} \in (x_a,x_b)$
and hence $\{\hat{g}=0\}$ contains $E_a \cap E_b$.\\ \\
\indent When $\{\hat{g}=0\}$ contains $E_a \cap E_b$, there exists
an affine chart in which one has\[\hat{g}=x_ag_a+x_bg_b+x_ax_bg_r,
\quad E_a \leftrightarrow x_a=0, \quad \quad E_b \leftrightarrow
x_b=0,\] with $g_a \notin (x_b)$ and $g_b \notin (x_a)$. If $m \in$
supp$(g)$ such that $\nu_a(m)=v_a$ and $\nu_b(m)=v_b$, Lemma $2$
implies that $\hat{m} \in$ supp$(x_ag_a) \cap$ supp$(x_bg_b)$ what
is
impossible. \vspace*{-0.5cm} \[\] \hfill $\blacksquare$\\ \\
\indent Now let $g \in \mathbb{P}F_{\underline{v}}$ and let $L$ be
the support of $g$. Taking the above characterization into account,
we see that if $\hat{g}$ does not contain $E_a \cap E_b$, then there
exists no $D \subset J$ for which $a \in D$, $b \notin D$ and
$\left(\cup_{d \in D}\Lambda_{L,d}\right) \cap \left(\cup_{d \notin
D}\Lambda_{L,d}\right)=\emptyset$. Note that the other implication
is false in general (see for example the constellation given above
with $a=1$ and $b=2$).
\\
We denote
\[Z_{\underline{v},L}:=\{g \in Z_{\underline{v}} \mid
\mbox{supp}(g)=L\} \mbox{ and } D_{\underline{v},L}:=\{\{\hat{g}=0\}
\cap \mathcal{D} \in D_{\underline{v}} \mid \mbox{supp}(g)=L\}.\]
Then we can write $Z_{\underline{v}}$ as a disjoint union
$\cup_{L}Z_{\underline{v},L}$ where for each $L$ holds that there is
no subset $D$ of $J$ satisfying condition $(2)$. Proposition $4$
tells us that the map
\begin{eqnarray*}
\psi_L: Z_{\underline{v},L} & \longrightarrow & D_{\underline{v},L}
\\
g & \longmapsto & \{\hat{g}=0\} \cap \mathcal{D}
\end{eqnarray*}
is a bijection and then part $1$ of the proof of Proposition $4$
allows us to conclude that the map
\begin{eqnarray*}
\phi: Z_{\underline{v}} & \longrightarrow & D_{\underline{v}}\\
g & \longmapsto & \{\hat{g}=0\} \cap \mathcal{D}
\end{eqnarray*}
is bijective.
\\ \\ \indent
Now by Lemma 2 one can see that $\phi$ is a homeomorphism. Then we
have that $\chi(Z_{\underline{v}})=\chi(D_{\underline{v}})$ which
completes the proof of Theorem $2$. \vspace*{-0.5cm} \[\] \hfill
$\blacksquare$
\\ \\
\noindent \textbf{Corollary.---} \emph{Let
$\underline{v}_1,\underline{v}_2, \cdots, \underline{v}_s$ be the
different columns of the linear system determined by the
constellation. Then}
\[P(\underline{t})=\frac{1}{(1-\underline{t}^{\underline{v}_1})^{\chi(D_{\underline{v}_1})}
\cdots
(1-\underline{t}^{\underline{v}_s})^{\chi(D_{\underline{v}_s})}}.\]
\emph{Proof.} When the cone $C$ associated to the linear system
$\mathcal{L}(\mathcal{C})$ is non-degenerated, then one has for each
$\underline{v}_i$ ($1 \leq i \leq s$) that
$\chi(D_{\underline{v}_i})=N(\underline{v}_i)=1$. If $C$ is
degenerated then we have for all but one $\underline{v}_i$ ($1 \leq
i \leq s$) that $\chi(D_{\underline{v}_i})=N(\underline{v}_i)=1$.
For the column $\underline{v}$ that appears more than once, one gets
$\chi(D_{\underline{v}})=N(\underline{v})=k$, with $k=d-s+1$.
\vspace*{-0.5cm} \[\] \hfill $\blacksquare$
\\ \\
As a consequence of Theorem $2$ and the corollary, we obtain that
the value $\chi(D_{\underline{v}})$ can be calculated from the
values $\chi(D_{\underline{v}_1}), \cdots,
\chi(D_{\underline{v}_s})$ and $k$.
\\ \\ \\ 

\footnotesize{

Ann Lemahieu, K.U.Leuven, Departement Wiskunde, Celestijnenlaan
200B, B-3001 Leuven, Belgium\\
\indent E-mail address: ann.lemahieu@wis.kuleuven.be

\begin{thebibliography}{DL3}

\bibitem[B,G-S]{BouvierGonzalezSprinberg} C. Bouvier and G. Gonzalez-Sprinberg, \emph{Syst\`{e}me g\'{e}n\'{e}rateur minimal,
diviseurs essentiels et G-d\'{e}singularisation de vari\'{e}t\'{e}s
toriques}, T\^{o}hoku Math. J. \textbf{47} (1995), 125-149.

\bibitem[B]{Bouvier} C. Bouvier, \emph{Diviseurs essentiels, composantes essentielles des vari\'{e}t\'{e}s
toriques singuli\`{e}res}, Duke Math. J. \textbf{91} (1998),
609-620.

\bibitem[C,D,G-Z1]{trio1} A. Campillo, F. Delgado, S.M. Gusein-Zade, \emph{The extended semigroup of a plane curve singularity},
Proc. Steklov Institue Math. \textbf{221} (1998), 139-156.

\bibitem[C,D,G-Z2]{trio2} A. Campillo, F. Delgado, S.M. Gusein-Zade, \emph{The Alexander polynomial
of a plane curve singularity via the ring of functions of it}, Duke
Math. J. \textbf{117} (2003), 125-156.

\bibitem[C,D,G-Z3]{trio3} A. Campillo, F. Delgado, S.M. Gusein-Zade, \emph{Integration with respect to the Euler characteristic over the
projectivization of the space of functions and the Alexander
polynomial of a plane curve singularity}, Russ. Math. Surv.
\textbf{55} (2000), 1148-1149.

\bibitem[C,D,G-Z4]{trio4} A. Campillo, F. Delgado, S.M. Gusein-Zade, \emph{Integrals with respect to the Euler characteristic
over the space of functions and the Alexander polynomial}, Proc.
Steklov Inst. Math. \textbf{238} (2002), 134-147.

\bibitem[C,D,G-Z5]{trio5} A. Campillo, F. Delgado, S.M. Gusein-Zade, \emph{The Alexander polynomial of a plane curve singularity and integrals
with respect to the Euler characteristic}, Internat. J. Math.
\textbf{14}(1) (2003), 47-52.

\bibitem[C,D,G-Z6]{trio6} A. Campillo, F. Delgado, S.M. Gusein-Zade, \emph{Poincar\'{e} series of a rational surface singularity}, Inv. Math. \textbf{155}
(2004), 41-53.

\bibitem[C,G-S,L-J]{clusters} A. Campillo, G. Gonzalez-Sprinberg, M. Lejeune-Jalabert, \emph{Clusters of infinitely near points}, Math. Ann. \textbf{306}
(1996), 169-194.

\bibitem[D,G-Z]{duo} F. Delgado, S.M. Gusein-Zade, \emph{Poincar\'{e} series for several plane divisorial valuations}, Proc.
Edinburgh Math. Soc. \textbf{46} (2003), 501-509.

\bibitem[D,L]{denefloeser} J. Denef, F. Loeser, \emph{Germs of arcs on singular algebraic varieties and motivic
integration}, Invent. Math. \textbf{135} (1999), 201-232.

\bibitem[I,K]{Ishii} S. Ishii, J. Koll\'{a}r, \emph{The Nash problem on arc families of singularities},
Duke Math. J. \textbf{120} (2003), no. 3, 601-620.

\bibitem[O]{Oda} T. Oda, \emph{Convex Bodies and Algebraic Geometry, An Introduction to the
Theory of Toric Varieties}, Springer-Verlag, Berlin, Heidelberg, New
York 1988.

\bibitem[V]{Viro} O.Y. Viro, \emph{Some integral calculus based on Euler characteristic}, Topology and Geometry - Rohlin seminar.
Lect. Notes Math. \textbf{1346}, 127-138. Berlin, Heidelberg, New
York: Springer 1988.

\end{thebibliography}
 \end{document}